\documentclass[reqno]{amsart}
\usepackage{graphics}
\numberwithin{equation}{section}

\renewcommand{\(}{\left(}
\renewcommand{\)}{\right)}

\begin{document}
\title{Verhulst's Logistic Curve}

\author{David M. Bradley}
\address{Department of Mathematics and Statistics\\
         University of Maine\\
         5752 Neville Hall\\
         Orono, Maine 04469--5752\\
         U.S.A.}
\email{bradley@math.umaine.edu, dbradley@member.ams.org}
\date{Submitted November 4, 1999.}



\maketitle

\section{Introduction}
Students tend to regard the elongated
``S-shaped''~\cite{HHG,Spieg,Strang} logistic curve of population
dynamics (fig.\ 1) as somewhat exotic.  It is typically derived by
applying the method of partial fractions to a separable differential
equation~\cite{Boy,CH,HHG,Rain,Spieg,Stew,Strang,Stroy,Will,Zill}.
My purpose here is to show how the logistic curve may be derived
more directly as a simple consequence of the more familiar
differential equation model for exponential decay, and that the
curve itself is nothing more than a familiar friend in disguise. The
disguise is removed by abandoning our fixation on the reference
point $(0,P_0)$, representing the initial population at time zero,
in favour of a much more natural choice. This illustrates an
important principle, namely that one should always adapt the
coordinates to the problem at hand.  In this case, a great deal is
simplified by relocating the origin more appropriately.

\section{Background}
Textbooks (see
eg.~\cite{Boy,CH,HHG,Rain,Stew,Strang,Stroy,Will,Zill}) typically
begin the discussion of population growth with the exponential
model
\begin{equation}
   \frac{1}{P}\frac{dP}{dt} = k \qquad\Longrightarrow
   \qquad P(t)=P_0e^{kt},\qquad
   t\ge 0,
\label{expgrow}
\end{equation}
in which the relative growth rate $k$ is a positive constant,
representing say the average birth rate.  Since unbounded growth is
unrealistic, more sophisticated models take into account limited
resources for reproduction. The logistic model, proposed by the
Belgian mathematical biologist Pierre F.\ Verhulst in
1838~\cite{Boy}, replaces the constant relative growth rate
in~(\ref{expgrow}) with a relative growth rate that decreases
linearly as a function of $P$:
\begin{equation}
   \frac{1}{P}\frac{dP}{dt} = k\(1-\frac{P}{M}\),\qquad k>0,\qquad
   0<P_0<M.
\label{logisticgrow}
\end{equation}
The constant $M$ represents the maximum sustainable population
beyond which $P$ cannot increase.  The dimensionless factor
$1-P/M$ in~(\ref{logisticgrow}) serves to diminish the relative
growth rate from $k$ down to zero as the population increases from
its initial level $P_0$ to $M$.

Although one can solve~(\ref{logisticgrow}) as a Bernoulli
differential equation by making the substitution
$P=1/y$~\cite{Pow}, for the most part texts
treat~(\ref{logisticgrow}) as a separable differential equation to
be solved by the method of partial
fractions~\cite{Boy,CH,HHG,Rain,Spieg,Stew,Strang,Stroy,Will,Zill}.
Either way, one obtains, after some algebraic simplifications, the
solution
\begin{equation}
   P(t)=\frac{M}{1+R_0e^{-kt}},\qquad\mbox{where}\qquad
   R_0=\frac{M-P_0}{P_0}.
\label{soln}
\end{equation}


\begin{figure}[ht]
   \begin{center}
   \scalebox{0.6}{\includegraphics{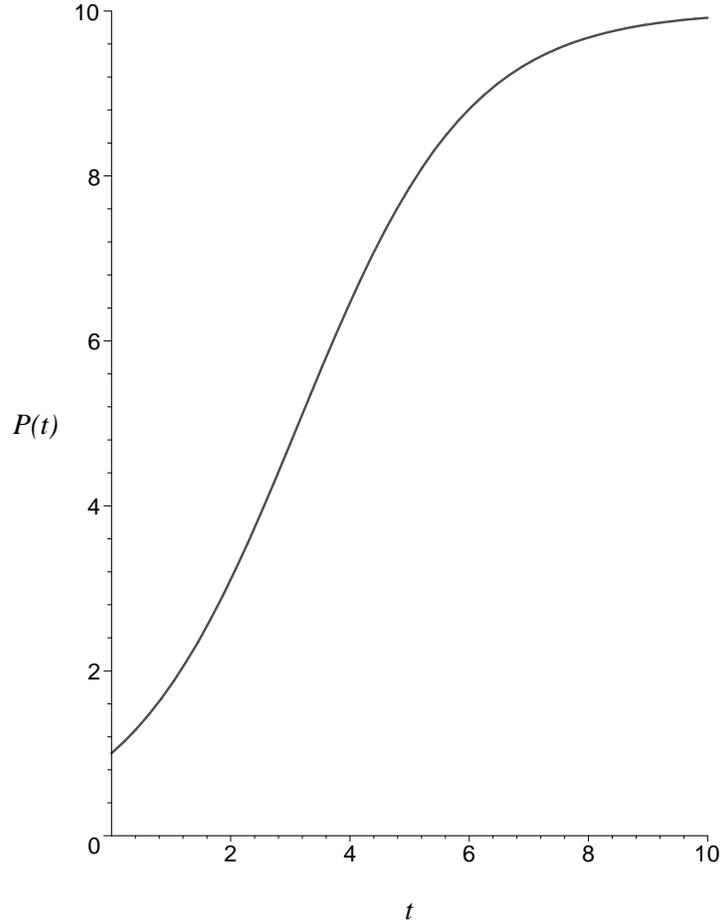}}
   \caption{Logistic Curve with $M=10$, $P_0=1$, $k=0.7$}
   \end{center}
\end{figure}

\section{Approach Via Exponential Decay}
Suppose instead of counting individuals, we count niches, viewing
$M$ as the maximum number of niches the ecosystem can support, and
$P$ as the number of niches currently occupied. Dimensional analysis
suggests that instead of studying $P$, we consider
\begin{equation}
   R=\frac{M-P}{P},
\label{ratio}
\end{equation}
the dimensionless ratio of available or vacant niches to niches
currently occupied.  From~(\ref{logisticgrow}), it readily follows
that
\begin{equation}
   \frac{1}{R}\frac{dR}{dt}=-k\qquad\Longrightarrow\qquad
   R(t)=R_0e^{-kt},
\label{Rsoln}
\end{equation}
where $R_0$ is as in~(\ref{soln}).  Expressing~(\ref{Rsoln}) in
terms of $P$, one arrives again at the solution~(\ref{soln}).

For the instructor who would like to discuss logistic growth but
would prefer to bypass separable differential equations and
partial fractions, it may be desirable to ``cut to the chase'' by
introducing the logistic model via~(\ref{ratio}) and~(\ref{Rsoln})
rather than via~(\ref{logisticgrow}).  In other words, although
the fact that $R$ satisfies the differential
equation~(\ref{Rsoln}) follows from its definition~(\ref{ratio})
and the logistic differential equation~(\ref{logisticgrow}), one
could equally well \emph{dispense with~(\ref{logisticgrow}) and
make~(\ref{Rsoln}) an assumption of the model}, thereby proceeding
more quickly to the solution.  Of course, there are pedagogical
advantages to either approach.  One aspect the approach
via~(\ref{ratio}) and~(\ref{Rsoln}) we are proposing has in its
favour is that logistic growth can be introduced in the standard
section on exponential growth and decay, with no loss in
continuity and without any additional background.

That $R$ decreases at a rate proportional to itself, i.e.\ satisfies
the differential equation~(\ref{Rsoln}), is intuitively plausible.
Initially we think of $P$ being much smaller than $M$, so that $R$
is much larger than $1$ and many niches are available relative to
the number currently occupied (a high niche vacancy rate).  We
should expect any species to take advantage of such a hospitable
climate for reproduction, and hence initially, $R$ should decrease
rapidly as $P$ increases.  However, as the number of vacancies
decreases, ($P$ gets close to $M$, $R$ gets close to zero) there are
relatively few available niches remaining.  In such an
\emph{inhospitable} climate, we should expect reproduction and hence
further growth to be difficult, and accordingly, $R$ should decrease
much more slowly. These considerations should be sufficient to
motivate the introduction of logistic growth via~(\ref{ratio})
and~(\ref{Rsoln}) to any calculus or differential equations class.

\section{Removing the Disguise}
From the viewpoint of an individual of the species attempting to
reproduce, one should expect a qualitative change in the
hospitality of the ecosystem near $R=1$, given the considerations
of the previous paragraph.  Motivated by these considerations, we
refer to an ecosystem as being {\it hospitable} or {\it
inhospitable} according to whether $R$ is greater or less than 1.
From~(\ref{Rsoln}), the transition from hospitable to inhospitable
occurs when
\begin{equation}
   t=\tau_0:=\frac1{k}\log R_0 =
   \frac1{k}\log\(\frac{M-P_0}{P_0}\),
   \qquad R=1,\qquad P=\tfrac12 M.
\label{taunaught}
\end{equation}
It is well-known that this is precisely the time at which $P$ is
increasing most rapidly, as can be seen by completing the square
in~(\ref{logisticgrow}):
\[
   \frac{dP}{dt}=\frac{Mk}{4}-\frac{k}{M}\(P-\frac{M}{2}\)^2.
\]
Because of the distinguished nature of the point $(\tau_0,M/2)$ it
seems more sensible to measure time from $\tau_0$ than from zero.
Certainly $t=0$ is completely arbitrary from the viewpoint of the
species, having more to do with whatever external forces (desire,
opportunity, availability of funding etc.) conspired to allow the
biologist or census taker to obtain an initial field count than
any essential features of the system.  Therefore, we consider
\[
   Q(\tau) := P(\tau_0+\tau) -\tfrac12 M,
\]
where $\tau$ measures time from $\tau_0$ and hence may be positive
or negative.  From~(\ref{soln}), we have
\[
   Q(\tau) = \frac{M}{1+R_0e^{-k(\tau_0+\tau)}}-\frac{M}2.
\]
Since $R_0e^{-k\tau_0}=1$, this simplifies to
\begin{equation}
   Q(\tau)=\frac{M}{1+e^{-k\tau}}-\frac{M}{2}
   =\frac{M}2\(\frac{1-e^{-k\tau}}{1+e^{-k\tau}}\)
   =\tfrac12{M}\tanh\(\tfrac12{k\tau}\),
\label{Qsoln}
\end{equation}
or in other words,
\begin{equation}
   P(t)=\tfrac12{M}\(1+\tanh\(\tfrac12 k(t-\tau_0)\)\),
\label{tanh}
\end{equation}
where $\tau_0$ is given by~(\ref{taunaught}). Thus, the mysterious
``S-shaped''~\cite{HHG,Spieg,Strang} logistic curve is nothing
more than a translate of our old and familiar friend, the
hyperbolic tangent.

\section{Addendum}

If $0\le P_0\le M$, then $\infty\ge R_0\ge 0$.  The boundary cases
$P_0=0$ and $P_0=M$ correspond to $R_0=\infty,P(t)\equiv0$ and
$R_0=0,P(t)\equiv M$, respectively.  To complete the analysis of
logistic growth, it is necessary to consider what happens when
$P_0$ lies outside the closed interval $[0,M]$, i.e.\ $R_0<0$. The
solution~(\ref{soln}) is valid for such $R_0$, but~(\ref{tanh})
was predicated on the assumption $R_0>0$ in the definition of
$\tau_0$.  Putting $S_0=-R_0$, we have from~(\ref{soln}) that
\[
   P(t)=\frac{M}{1-S_0e^{-kt}},\qquad S_0>0.
\]
In this case, we define
\begin{equation}
   \tau_0 := \frac{1}{k}\log S_0
   \qquad\mbox{so that}\qquad e^{k\tau_0}=S_0=\frac{P_0-M}{P_0}.
\label{newtau}
\end{equation}
A calculation analogous to~(\ref{Qsoln}) reveals that
\[
   H(\tau) := P(\tau_0+\tau)-\tfrac12 M
   =\tfrac12 M\coth\(\tfrac12 k\tau\),
\]
or
\begin{equation}
   P(t) = \tfrac12 M\(1+\coth\(\tfrac12 k(t-\tau_0)\)\),
\label{coth}
\end{equation}
where now $\tau_0$ is given by~(\ref{newtau}).

\begin{figure}[hb]
   \begin{center}
   \scalebox{0.6}{\includegraphics{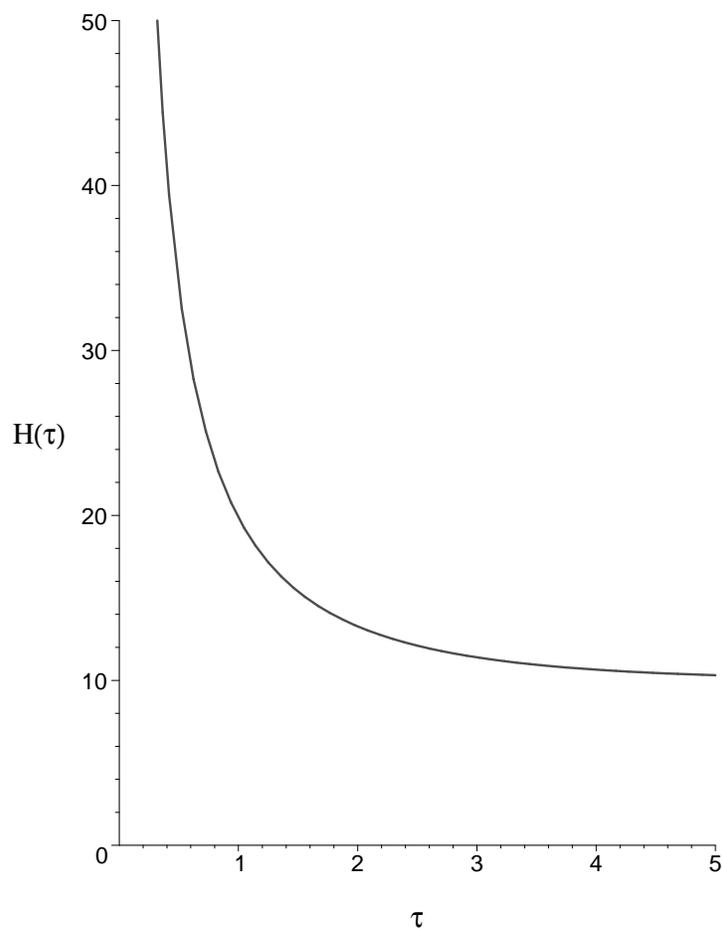}}
   \caption{Logistic Curve with $M=10$, $P_0=30$, $k=0.7$}
   \end{center}
\end{figure}

If $P_0>M$, then $-1\le R_0<0$, $0<S_0=-R_0\le 1$, and $-\infty
<\tau_0\le 0$. Therefore, for $t\ge 0$ we are on the upper arch of
the hyperbolic cotangent, with population decreasing exponentially
to $M$ as $t\to\infty$ (fig.\ 2).  In the less biologically
meaningful case $P_0<0$, we have $-1\ge R_0>-\infty$, $1\le
S_0=-R_0<\infty$ and $0\le \tau_0<\infty$.  As $t$ increases from
zero to $\tau_0$, the rightmost portion of the lower arch of the
hyperbolic cotangent is traversed, sending the population to minus
infinity.  The asymptote is then crossed and we skip over to the
upper arch, the population reverting to its behaviour in the
previous case.

\end{document}